\documentclass[11pt,reqno]{amsart}

\usepackage{amsfonts}
\usepackage{amsmath}
\usepackage{amssymb, esint}
\usepackage{amsthm,color}
\usepackage{enumerate}
\usepackage{enumitem}
\usepackage{mathtools}
\usepackage{url}
\usepackage[hidelinks, bookmarksdepth=3]{hyperref}
\usepackage{xcolor}
\usepackage{cancel}
\usepackage{ulem}
\usepackage{caption,subcaption}
\usepackage{tikz}
\usepackage{tikz-3dplot}

\usepackage{mathrsfs}
\usepackage{bm}

\usepackage{longtable}

\usepackage{soul}

\tdplotsetmaincoords{70}{110}

\usepackage{geometry}
\newcommand{\R}{{\mathbb R}}

\newcommand{\N}{{\mathbb N}}

\newcommand{\D}{{\mathbb D}}
\newcommand{\cE}{{\mathcal E}}

\newcommand{\cN}{{\mathcal N}}

\newcommand{\cA}{{\mathcal A}}

\newcommand{\cH}{{\mathcal H}}

\newcommand{\cF}{{\mathcal F}}

\def\0{{\mathbf 0}}
\def\loc{{\textup{loc}}}
\def\id{{\rm id}}

\def \mb{\mathbb}

\newcommand{\e}{\varepsilon}
\newcommand{\vp}{\varphi}

\newcommand{\supp}{\operatorname{supp}}

\newcommand{\dist}{\operatorname{dist}}

\newcommand{\rank}{\operatorname{rank}}

\theoremstyle{plain}
\newtheorem{thm}{Theorem}[section]

\newtheorem{cor}[thm]{Corollary}
\newtheorem{lem}[thm]{Lemma}
\newtheorem{prop}[thm]{Proposition}
\newtheorem{prob}[thm]{Problem}

\newcommand{\thistheoremnames}{}
\newtheorem*{genericthms}{\thistheoremnames}
\newenvironment{para*}[1]
  {\renewcommand{\thistheoremnames}{#1}%
   \begin{genericthms}}
  {\end{genericthms}}

\newtheorem{defn}[thm]{Definition}

\theoremstyle{remark}

\newtheorem{ex}[thm]{Example}

\newtheorem{claim}{Claim}
\newtheorem*{claim*}{Claim}

\newtheorem{rem}[thm]{Remark}

\numberwithin{equation}{section}

\title[Constraint Maps]{Constraint Maps: Insights and Related Themes}

\author[A.\ Figalli]{Alessio Figalli}
\address{Department of Mathematics, ETH Z\"urich,  Raemistrasse 101, 8092 Z\"urich, Switzerland }\email[A.\ Figalli]{alessio.figalli@math.ethz.ch}

\author[A.\ Guerra]{Andr\'e Guerra}
\address{Department of Pure Mathematics and Mathematical Statistics,  University of Cambridge,  Wilberforce Rd, Cambridge CB3 0WB, UK}\email[A.\ Guerra]{adblg2@cam.ac.uk}

\author[S.\ Kim]{Sunghan Kim}
\address{Department of Mathematics, Uppsala University, S-751 06 Uppsala, Sweden}\email[S.\ Kim]{sunghan.kim@math.uu.se}

\author[H.\ Shahgholian]{Henrik Shahgholian}
\address{Department of Mathematics, KTH Royal Institute of Technology, 100 44 Stockholm, Sweden } \email[H.\ Shahgholian]{henriksh@kth.se}

\begin{document}

\maketitle

\begin{abstract}
This paper explores recent progress related to constraint maps. Building on the exposition in \cite{FGKS3}, our goal is to provide  a clear and accessible account of some of the more intricate arguments behind the main results in this work. Along the way, we include several new results of independent value. In particular, we give optimal geometric conditions on the target manifold that guarantee a unique continuation result for the projected image map. We also prove that the gradient of a minimizing harmonic map (or, more generally, of a minimizing constraint map) is an $A_\infty$-weight, and therefore satisfies a strong form of the unique continuation principle. In addition, we outline possible directions for future research and highlight several open problems that may interest researchers working on free boundary problems and harmonic maps.
\end{abstract}

\medskip

\begin{center}
    \textit{
   To Sandro Salsa, on his 75th birthday,
   in recognition of his profound\\ contributions  to the theory and development of free boundary problems.
   }
\end{center}

\setcounter{tocdepth}{2}

\tableofcontents

\section{Introduction}

This article reviews recent developments on {\it Constraint Maps}. Extending the exposition in \cite{FGKS3}, we discuss several of the more intricate arguments that underlie the results in \cite{FKS, FGKS, FGKS2} and present a few new theorems that we believe are of independent interest. We conclude by outlining plausible directions for future work and listing open questions that are likely to interest researchers in free-boundary problems and harmonic maps.

A \emph{constraint map} is a critical point of the Dirichlet energy
\[
\mathcal{E}[u]=\int_{\Omega}|Du|^{2}\,dx,
\]
subject to the image constraint \(u\in W^{1,2}(\Omega;\overline{M})\).  
Here \(\Omega\subset\mathbb{R}^{n}\) (\(n\ge 1\)) is a bounded open set and \(M\subset\mathbb{R}^{m}\) (\(m\ge 2\)) is a smooth domain whose complement \(M^{c}\) is compact and connected; following the classical obstacle problem~\cite{Caf98}, we call \(M^{c}\) the \emph{obstacle}.  
We concentrate on energy--minimising maps with fixed boundary data, since the problem is already both rich and delicate for this class.

Existence of minimisers with prescribed boundary values follows from the direct method of the calculus of variations whenever the admissible class is non-empty---for instance, when the boundary data take values in \(\overline{M}\) and \(n\ge 3\).\footnote{See Remark~2.9 in~\cite{FGKS2} for further discussion.}  
Uniqueness, however, generally fails; explicit counterexamples appear even in one dimension (see~\cite{LL} for a detailed analysis).  
Throughout, it is helpful to view constraint maps as parametrizations of a given geometric object under constraints, a perspective that clarifies both the similarities and the differences with the scalar obstacle problem (cf. Section 2.5 of \cite {FGKS2}).

The subtle interplay between free boundaries and mapping singularities has been explored in recent work~\cite{FGKS,FGKS2}, concentrating on two principal types. 
The first consists of \textit{discontinuity points}---singularities that arise from a topological mismatch between the domain and the target and therefore cannot be removed by small perturbations.  Such defects lie at the core of harmonic-map theory.  
The second involves \textit{branch points}, which may signal either a genuine singularity in the map’s image or a singular parametrisation of an otherwise regular image. This phenomenon is reminiscent of the classical theory of minimal surfaces.

In order to make the last connection more explicit, let us briefly discuss constraint maps from the perspective of minimal surface theory, as it played an important historical role in the subject. The main observation is that \textit{weakly conformal} constraint maps can be considered as parametric minimal surfaces with obstacles. To be more precise, let us recall the setting of the classical Plateau problem. Let $\Gamma$ be a Jordan curve in $\R^m$ ($m\geq 3$). We aim at finding the area-minimizing surface spanning $\Gamma$. The Plateau problem was solved in the early 1930s independently by Douglas and Rad\'o, see e.g.\ \cite{Colding2011} for a modern exposition. The solutions are given by considering energy-minimizing parametrizations $u$ from the unit disk $\D \subset \R^2$ to $\R^m$ that are monotone on the boundary onto $\Gamma$ and satisfy the so-called three point condition\footnote{That is, we prescribe three points on the boundary of $\partial\D$, say $(-1,0)$, $(0,1)$, and $(1,0)$, and also three points on the image $\Gamma$, say $\xi_1$, $\xi_2$, and $\xi_3$, and ask $u(-1,0) = \xi_1$, $u(0,1) = \xi_2$, and $u(1,0) = \xi_3$.}. The energy-minimality under the boundary condition ensures that the parametrizations are weakly conformal (i.e., $|u_1| = |u_2|$ and $u_1\cdot u_2 = 0$ a.e., where $u_i = \frac{\partial u}{\partial x_i}$, $i=1,2$). Then the (generalized) surface area $\cA$ of the image $\Sigma = u(\D)$ equals its energy (see \cite[Page 135]{Colding2011}): 
$$
\cA(\Sigma)= \int_\D \sqrt{|u_1|^2|u_2|^2 - (u_1\cdot u_2)^2}\,dx = \frac{1}{2} \cE(u,\D).
$$
Now let us place an obstacle $M^c$ within $\Gamma$ (or, more precisely, consider $\Gamma$ that do not intersect $M^c$). Tomi \cite{T1,T2} and Hilderbrandt \cite{H1,H2}  independently in the early 1970s studied the above Plateau problem for surfaces constrained to lie outside of ${\rm int}(M^c)$; in particular, the surface can lie on the boundary of the obstacle. The solutions are, by our definition above, energy minimizing weakly conformal constraint maps, and they also yield parametric minimal surfaces with an obstacle. We also note that minimal surfaces with obstacles are a classical subject, having been investigated both in the non-parametric setting \cite{Giusti2010,Kinderlehrer1973} and in the setting of sets of finite perimeter \cite{Barozzi1982,Miranda1971,Tamanini1982}.

In the 1980s, the theory of Tomi and Hildebrandt was then generalized by Duzaar and Fuchs in a series of works, see e.g.\ \cite{D,DF}. There are few new results in the topic after the 1990s, until  recently the authors revitalized the problem in a series of works \cite{FKS, FGKS, FGKS2}. For a more detailed historical recap of the development of this theory we refer the reader to the introduction of \cite{FGKS2} and the  expository article \cite{FGKS3}. 

\bigskip

The present note discusses the unique continuation principle (UCP) for constraint maps and highlights the relevance of these properties to recent results in~\cite{FGKS3}, while also deriving several partial results of independent interest. 
As a consequence of the $\e$-regularity theory for minimizing constraint maps and the classical unique continuation principle, one  concludes the following theorem.

\begin{thm}[Classical UCP \cite{A}]\label{thm:classic}
    Let $u\in W^{1,2}(\Omega;\overline M)$ be a minimizing constraint map. If $|Du|=0$ on an open subset of $\Omega$  then $u$ is constant in $\Omega.$
\end{thm}

In this paper, we discuss finer UCP for constraint maps. Since $\partial M$ is smooth, there is a tubular neighborhood $\cN(\partial M)$ of $\partial M$ where we can define the nearest point projection $\Pi$ to $\partial M$, which is then smooth in $\cN(\partial M)$. A basic question is whether the projected image $\Pi\circ u$ has the UCP, when $u$ is a minimizing constraint map. The following gives a rather complete answer to this question:

\begin{thm}\label{thm:UCPPi}
Let $u\in W^{1,2}(\Omega;\overline M)$ be a minimizing constraint map. Suppose that there is a ball $B\subset \Omega$ such that  $u\in\cN(\partial M)$ a.e.\ in $B$ and that $|D(\Pi\circ u)| = 0$ on a set $F\subset B$ of positive measure. Denoting by $\nu$ the normal to $\partial M$ pointing towards $M$, if
\begin{equation}
    \label{eq:geometriccondintro}
    \{y + t\nu(y): t> 0\}\cap M^c =\emptyset \quad \text{ for every } y\in\partial M,
\end{equation}
 then $\Pi\circ u$ is constant in $\Omega$.  Furthermore, the conclusion fails in the absence of condition \eqref{eq:geometriccondintro}.
\end{thm}

\begin{rem}
    The condition \eqref{eq:geometriccondintro} holds for instance if our obstacle $M^c$ is convex, but it also works for certain non-convex obstacles.  
\end{rem}

There are many other basic questions concerning $\Pi\circ u$ which are not yet solved; for instance, it is not known what the optimal regularity of this map is, see Problem \ref{prob:Pi} below and \cite{FKS}.

\medskip

Our second main result concerns instead \textit{quantitative} versions of the UCP. A natural way of formulating such a result involves the concept of $A_\infty$-weight, cf.\ Definition \ref{def:Ainfty} below. In a nutshell, an $A_\infty$-weight is a non-negative function that satisfies a reverse H\"older inequality with non-increasing supports. The role of $A_\infty$-weights in proofs of UCP results for elliptic PDE is classical, see e.g.\ \cite{Garofalo1986} for the case of linear elliptic equations with Lipschitz coefficients. In contrast, when studying constraint maps one is forced to consider systems of equations with \textit{measurable} coefficients, since the term $A_u(Du,Du)$ appearing in the Euler--Lagrange system \eqref{eq:EL} is not even bounded around the singular set of $u$. 

We are aware of only two results concerning UCP for solutions of elliptic PDE with measurable coefficients. The first such result is the weak Harnack inequality, which holds for \textit{non-negative} solutions of \textit{scalar} elliptic equations both in divergence and non-divergence form, respectively due to the works of Moser (after De Giorgi and Nash) and Krylov--Safonov, see e.g.\ \cite{Kenig1993,Moser1961}. We note that the non-negativity assumption is crucial, as the UCP fails already for sign-changing solutions of elliptic equations with $\alpha$-H\"older coefficients for any $\alpha<1$, see \cite{Plis1963}. The second such result is due to Gehring \cite{Gehring1973a} in the context of quasiregular maps, which are solutions of certain first-order elliptic systems, see e.g.\ \cite{Iwaniec2001} for a modern PDE-based exposition. In fact, Gehring's work pioneered the study of reverse H\"older inequalities and $A_\infty$-weights.

Here we prove that minimizing harmonic maps, and more generally minimizing constraint maps, also give rise to $A_\infty$-weights (see Definition \ref{def:Ainfty}):

\begin{thm}\label{thm:HMs}
    Let $N\subset \R^m$ be a smooth compact manifold without boundary, and let $u\in W^{1,2}(\Omega,N)$ be a minimizing harmonic map into $N$. Then $|Du|$ is an $A_\infty$-weight.
\end{thm}

\begin{rem}\label{rem:manifoldswithbdry}
    The same theorem holds if $N\subset \R^m$ is a sufficiently ``well-behaved'' complete manifold \textit{with boundary}, e.g.\ if $N=\overline M$ is the complement of a smooth bounded domain, as in this paper. In general one needs to ask a condition on $N$ similarly to Luckhaus' paper \cite{L}, i.e.\ one needs that
    $$\limsup_{r\to 0}\left\{\frac{|x'-y'|}{|x-y|}: \dist(x,N)+\dist(y,N)<r, x'\in \Pi_N(x),y'\in \Pi_N(y)\right\}= 1,$$
    where $\Pi_N(z)=\{z'\in N:|z-z'|=\dist(z,N)\}$ is the set of nearest point projections of $z$ in $N$.
\end{rem}

Theorem \ref{thm:HMs} and the previous remark give a strengthening of one of the main results in \cite{FGKS2}, in particular Theorem 1.2 therein, and here we give a more streamlined proof of that result using the theory of $A_\infty$-weights.

\medskip

We conclude this introduction by describing the organization of the paper. In Section \ref{sec:prelims} we gather some introductory material in order to fix the notation and basic properties of constraint maps. Then, in Section \ref{sec:radil}, we discuss the minimizing properties of some concrete examples of constraint maps with symmetries. In the rest of the paper we then discuss the UCP in detail, beginning in Section \ref{sec:proj} with the UCP for the projected image and the proof of Theorem \ref{thm:UCPPi}. In Section \ref{sec:ucpfb} we explain how the UCP is related to regularity of constraint maps around their free boundaries, and in Section \ref{sec:Ainfty} we prove the quantitative form of UCP asserted in Theorem \ref{thm:HMs}. Finally, in Section \ref{sec:probs} we provide a list of open problems, hoping to invite more people to study this interesting subject.

\section{Notation and terminology}
\label{sec:prelims}

 We mostly follow  the notation and terminologies  used in  \cite{FGKS2}. In particular, we denote by $u^{-1}(M)$ the non-coincidence set, by $u^{-1}(\partial M)$ the coincidence set, and the interface between these two sets, namely $\partial u^{-1}(M)$, as the free boundary. Since $\partial M$ will always be assumed to be smooth, there is a tubular neighborhood $\cN(\partial M)$ of $\partial M$ where we can define the signed distance function $\rho$ to $\partial M$ (which takes positive values in $M$) and the nearest point projection $\Pi$ to $\partial M$: on $\cN(\partial M)$, both $\rho$ and $\Pi$ are smooth. We also denote the unit normal to $\partial M$ by $\nu$ (pointing towards $M$) and by $A$ the second fundamental form of $\partial M$.

 Throughout most of our discussion, $u$ is going to be an \textit{energy minimizing} constraint map. However, it is also natural to consider larger classes of constraint maps. The largest class of maps that we can reasonably consider is the class of \textit{weakly constraint maps}
\cite{CM}: this is the family of maps that are critical points of the Dirichlet energy in the following sense:
 $$\liminf_{v\to u, v \in W^{1,2}(\Omega;\overline M)}\frac{\cE[u]-\cE[v]}{\|u-v\|_{W^{1,2}(\Omega)}}\geq 0.$$
Under reasonable conditions on $M$ (which hold for the class of targets we consider here), one can show that weakly constraint maps are exactly the solutions of the system
$$\Delta u = \nu_u \,\lambda,$$
where $\lambda$ is a Radon measure such that
\begin{equation}
\label{eq:lambda}
0 \leq \lambda\leq\left( -A_u(Du,Du)\cdot \nu_u\, \right) \chi_{u^{-1}(\partial M)},
\end{equation}
cf.\ \cite[Theorem 3.1]{CM}. Whenever $u\in C(\Omega\setminus \Sigma)$ where $\Sigma$ is a null-set, one can verify that equality holds in the second inequality in \eqref{eq:lambda}.

The class of weakly constraint maps certainly contains the class of weakly harmonic maps into $\partial M$: these are simply those constraint maps for which $u^{-1}(M) =\emptyset$. However, as soon as $n\geq 3$, and even when $\partial M=\mb S^2$, there is no partial regularity theory for such maps \cite{Riviere1995}. Thus one is forced to work with more restricted classes of maps. The results of this paper all hold for the following class of maps:
\begin{defn}\label{def:F}[Tangent maps and the class $\cF(\Omega;\overline M)$.]

    Let $u \in W^{1,2}(\Omega;\overline M)$ be given, and define $u_{x_0,r}(x) := u(x_0 + rx)$ whenever $x_0 + rx\in\Omega$.
    \begin{itemize}
        \item Given $x_0\in \Omega$, we call $\vp$ a tangent map of $u$ at $x_0$ if $u_{x_0,r_k} \rightharpoonup \vp$ in $W_\loc^{1,2}(\R^n;\R^m)$ along a sequence $r_k\downarrow 0$. 
        \item We say that $u\in \cF(\Omega;\overline M)$ if the following hold:
\begin{enumerate}
    \item[(i)] $u$ is a weak solution to 
    \begin{equation}
        \label{eq:EL}
        \Delta u = A_u(Du,Du)\chi_{u^{-1}(\partial M)}\quad\text{in }\Omega,
    \end{equation}
    where $A_{u(x)}$ is the second fundamental form of $\partial M$ at $u(x)$. 
    \item[(ii)] Given $x_0\in\Omega$ and $0<r< R< \dist(x_0,\partial\Omega)$,
    $$
    E(u,x_0,r):=r^{2-n}\int_{B_r(x_0)} |Du|^2\,dx \leq R^{2-n}\int_{B_R(x_0)} |Du|^2\,dx. 
    $$
    \item[(iii)] $u$ is continuous at $x_0$, i.e.\ $x_0\not \in\Sigma(u)$,  if and only if
    $$
    \lim_{r\downarrow 0} r^{2-n}\int_{B_r(x_0)} |Du|^2\,dx = 0. 
    $$
    \item[(iv)] If $\vp$ is a tangent map of $u$ at $x_0 \in \Omega$,  then there is a sequence $r_k\downarrow 0$ such that $u_{x_0,r_k} \to \vp$ in $W_\loc^{1,2}(\R^n;\R^m)$; moreover, $\vp\in \cF(B_R;\overline M)$ for all $R >0$.
\end{enumerate} 
    \end{itemize}
\end{defn}

\begin{rem}
    As we noted before, (i) holds for any weakly constraint map  whose discontinuity set is of measure zero. Condition (ii) is simply the monotonicity formula and it is satisfied for maps which are critical points for the Dirichlet energy with respect to inner-variations. Now (iii) is a qualitative version of the $\e$-regularity theorem. Stationary constraint maps (i.e.\ maps which are critical points for both inner and outer variations) satisfy (i)--(iii); the $\e$-regularity theorem is e.g.\ a consequence of \cite{Riviere2008}. However, the compactness assertion in (iv) is in general false for stationary maps. 
    \end{rem}

By \cite{D,L}, minimizing constraint maps satisfy conditions (i)--(iv). In fact, since these conditions are local, the same is true for locally minimizing constraint maps:\\
We say that $u\in W^{1,2}(\Omega;\overline M)$ is \textit{locally minimizing}, if for every $x\in\Omega$ there is a ball $B\subset\Omega$ centered at $x$ such that $\int_B|Du|^2\,dx\leq \int_B |Dv|^2\,dx$ for all $v\in W_u^{1,2}(B;\overline M)$. We will come back to the distinction between globally and locally minimizing maps in the next section. 

The above discussion leads us naturally to the following interesting question:\\
\textit{Is there a class of stationary constraint maps, strictly containing the class of locally minimizing maps, for which conditions (i)--(iv) hold?}\\
The theory of harmonic maps, see e.g.\ \cite[\S 3.4]{LW}, suggests that the class of \textit{stable-stationary maps} could be such a class, provided the target satisfies suitable geometric conditions, which hold e.g.\ when the target is $\mb S^{m-1}$ with $m\geq 4$. For harmonic maps into $\mb S^{m-1}$, the definition of stability is clear: one simply requires that
    $$0\leq \frac{d^2}{d t^2}\cE\bigg(\frac{u+t\phi}{|u+t\phi|}\bigg)\bigg|_{t=0} \quad \forall \phi \in C^\infty_c(\Omega,\R^m).$$
However, even when the obstacle is $\mb B^m$, constraint maps are in general at most $C^{1,1}$, and in particular they are not $C^2$. Thus a novel definition of stability for such maps is required; see Problem 7.5 in Section \ref{sec:probs}.

\section{On the minimality of radial and equatorial maps}\label{sec:radil}

In \cite[\S 2.4]{FGKS2} we provided some examples of radial constraint maps that are locally minimizing.  However, the most interesting question in this regard is whether those maps are actually \textit{globally} minimizing. Here we prove that this is indeed the case in large dimensions, i.e.\ whenever $n\geq 7$. This dimensional restriction is related to the sharp constant in Hardy's inequality, which plays an important role in the minimality of certain harmonic maps into spheres, cf.\ the classical work \cite{JK} as well as \cite{Hong2000} for a simpler approach, which we follow here.

\begin{prop}
\label{prop:radial}
    Let $n\geq 7$, $a\in(0,1)$, and let $u\in W_\id^{1,2}(B_1;B_a^c)$ be the radial weak solution to 
    \begin{equation}
        \label{eq:radial}
    \Delta u = - \frac{|Du|^2}{a^2}u \chi_{\{|u| = a\}}\quad\text{in }B_1.
    \end{equation}
    Then $u$ is the unique minimizing map in $W_\id^{1,2}(B_1;B_a^c)$.
\end{prop}

\begin{proof}
    Let us recall the sharp Hardy inequality: for every $\phi\in W_0^{1,2}\cap L^\infty(B_1)$, 
    \begin{equation}
        \label{eq:hardy}
        \frac{(n-2)^2}{4}\int_{B_1} \frac{\phi^2}{|x|^2}\,dx \leq \int_{B_1} \bigg|\frac{\partial\phi}{\partial r}\bigg|^2\,dx,
    \end{equation}
    where $\frac{\partial}{\partial r}$ is the directional derivative in the direction $|x|^{-1}x$, and equality holds if and only if $\phi=0$. Note that, by \cite[Example 2.10]{FGKS2}, radial weak solutions to \eqref{eq:radial} are unique\footnote{The radial profile solves a 2nd order ODE; see \cite{FGKS2} for more details.} and they take the form
    $$u(x) = w(|x|)\frac{x}{|x|}, \qquad w(r) = \begin{cases}
        a & \text{if }  r\in [0,r_a]\\
        t_a r + (1-t_a)r^{1-n} & \text{if } r\in (r_a,1]
    \end{cases},$$
    where the parameters $t_a,r_a\in(0,1)$ can be uniquely determined by the overdetermined condition on the free boundary. Note in particular that $\overline{B_{r_a}} = \{|u|=a\}$ and
    \begin{equation}
        \label{eq:radial-re}
    |Du(x)|^2 = \frac{(n-1)a^2}{|x|^2} \quad\text{in }B_{r_a}\setminus \{0\}.
    \end{equation}
    Now choose an arbitrary competitor $v\in W_\id^{1,2}(B_1;B_a^c)$. Then, thanks to the constraint $\min\{|u|,|v|\}\geq a$ a.e.\ in $B_1$, we obtain 
    \begin{equation}
        \label{eq:RHS}
\int_{B_1} \frac{|u-v|^2}{|x|^2}\,dx \geq 2 \int_{B_{r_a}} \frac{a^2 - u\cdot v}{|x|^2}\,dx =: I.
    \end{equation}
    Note that $I \neq 0$ implies $u\not\equiv v$ in $B_1$, as $\overline{B_{r_a}} = \{|u|=a\}$. Also, thanks to \eqref{eq:radial} and \eqref{eq:radial-re}, it follows that
    \begin{equation}
        \label{eq:LHS}
        \begin{aligned}
        \int_{B_1} |Du|^2 - |Dv|^2 + |D(u-v)|^2\,dx &= 2\int_{B_1} Du : D(u-v) \,dx \\
        & =  \frac{2}{a^2}\int_{B_{r_a}} |Du|^2 u\cdot (u-v)\,dx \\
        & =  2(n-1)\int_{B_{r_a}} \frac{u\cdot (u-v)}{|x|^2}\,dx= (n-1)I.
        \end{aligned}
    \end{equation}
    Now if $I \leq 0$, then \eqref{eq:LHS} directly yields that 
    $$
    \int_{B_1} |Du|^2- |Dv|^2\,dx \leq -\int_{B_1} |D(u-v)|^2\,dx \leq 0,
    $$
    so, $u$ is the unique minimizer within the subclass of functions $v \in W_\id^{1,2}(B_1;B_a^c)$ for which $I\leq 0$.
    
    In contrast, if $I > 0$ (and hence $u\not\equiv v$ in $B_1$), then combining \eqref{eq:LHS} with \eqref{eq:RHS}, and utilizing the Hardy inequality \eqref{eq:hardy} with $\phi = u^i - v^i \in W_0^{1,2}\cap L^\infty(B_1)$ for every $i=1,2\cdots,n$, we see that 
    $$
    \int_{B_1} |Du|^2 - |Dv|^2\,dx = (n-1)I-\int_{B_1} |D(u-v)|^2\,dx \leq \bigg(n-1 - \frac{(n-2)^2}{4}\bigg) I < 0,
    $$
    where the last inequality follows from $n \geq 7$ and $ I < 0$. (Note that this is the only place where the dimension comes into play.) Thus, the energy of $u$ is also strictly minimal among all competitors $c \in W_\id^{1,2}(B_1;B_a^c)$ satisfying $I > 0$.

    Combining the above observations, we conclude that $u$ is the unique global minimizer within the entire class $W_\id^{1,2}(B_1;B_a^c)$, provided that $n \geq 7$.
\end{proof}

We next consider the minimality of equator maps with respect to ellipsoidal obstacles, see \cite{Baldes1984,Helein1988} for related results in the harmonic map setting.

\begin{prop}
    \label{prop:equator}
    Fix $n\geq 7$ and $a\in(0,1)$. Let $u\in W_\id^{1,2}(B_1;B_a^c)$ be the radial weak solution to \eqref{eq:radial}. Then the map $(u,0)$ is also a unique minimizing map in $W_{(\id,0)}^{1,2}(B_1; E_{a,\lambda}^c)$, where $E_{a,\lambda}$ is the $(n+1)$-dimensional ellipsoid given by $\{(y,\xi): |y|^2 + \lambda^{-2}|\xi|^2 < a^2\}$, provided $\lambda^2 \geq \frac{4(n-1)}{(n-2)^2}$.
\end{prop}

\begin{proof}
    The proof here is essentially the same as that of Proposition \ref{prop:radial}. Let $V = (v,\xi)\in W_{(\id,0)}^{1,2}(B_1;E_{a,\lambda}^c)$ with $v = (v^1,\cdots,v^n)$, and write $U = (u,0)$. Then, instead of \eqref{eq:RHS}, we have
    $$
    \int_{B_1}\frac{|U-V|^2}{|x|^2}\,dx \geq \int_{B_{r_a}} \frac{a^2 - 2u\cdot v + |V|^2}{|x|^2}\,dx =: II.
    $$
    If $\lambda \geq 1$, then it follows from the constraint $|V|^2 \geq |v|^2 + \lambda^{-2}|\xi|^2 \geq a^2$ that 
    $$
    II \geq 2 \int_{B_{r_a}} \frac{a^2 - u\cdot v}{|x|^2}\,dx = I,
    $$ 
    where $I$ is as in \eqref{eq:RHS}. On the other hand, if $0<\lambda \leq 1$, we have $|V|^2  \geq \lambda^2|v|^2 + \xi^2 \geq \lambda^2 a^2$, whence
    $$
    II \geq 2\lambda^2 \int_{B_{r_a}} \frac{a^2 - u\cdot v}{|x|^2}\,dx = \lambda^2 I. 
    $$
    Note that the computation in \eqref{eq:LHS} remains the same for $U,V$ in place of $u,v$, since it only involves the system \eqref{eq:radial} that $u$ solves:
    $$
    \int_{B_1} |DU|^2 - |DV|^2 + |D(U-V)|^2\,dx = (n-1) I. 
    $$
    The rest of the proof is now the same: more precisely, when applying Hardy inequality, we use the restriction $n\geq 7$ when $\lambda \geq 1$, and the restriction $\lambda^2 \geq \frac{4(n-1)}{(n-2)^2}$ when $0<\lambda\leq 1$. We omit the details.  
\end{proof}

\section{Unique continuation property of the projected image}\label{sec:proj}

In this section we study the UCP of the projection map $\Pi\circ u$ when $u$ is a constraint map. In particular, we will prove Theorem \ref{thm:UCPPi}. 

Recall that $\Pi$ is the nearest point projection onto $\partial M$, which is well-defined and smooth in a tubular neighborhood $\cN(\partial M)$ of $\partial M$. 
We begin with the following more general version of the positive statement in Theorem \ref{thm:UCPPi}.

\begin{prop}\label{prop:proj-ucp}
 Let $u\in W^{1,2}(\Omega;\overline M)$ be a weakly constraint map such that $u\in C(\Omega\setminus\Sigma)$ for a closed null set $\Sigma$ that does not disconnect any open-connected set in $\Omega$. Suppose that there is a ball $B\subset \Omega$ such that  $u\in\cN(\partial M)$ a.e.\ in $B$ and that $|D(\Pi\circ u)| = 0$ on a set $F\subset B$ of positive measure. If $M$ satisfies
\begin{equation}
    \label{eq:geometriccond}
    \{y + t\nu(y): t> 0\}\cap M^c =\emptyset \quad \text{ for every } y\in\partial M,
\end{equation}
 then $u = y_0 + h\,\nu(y_0)$ in $\Omega$ for some point $y_0\in\partial M$ and a nonnegative harmonic function $h$ in $\Omega$.  In particular either $h>0$ or $h\equiv 0$.
\end{prop}

Note that the assumptions on $\Sigma$ are satisfied in particular whenever $\dim_{\cH}(\Sigma)<n-1$, see e.g.\ \cite[Theorem IV.4 and p.\ 107]{Hurewicz1942}.

\begin{proof}
By our assumptions and the discussion in Section \ref{sec:prelims}, $u$ satisfies
\begin{equation}
    \label{eq:u-pde}
    \Delta u = A_u(Du,Du)\chi_{u^{-1}(\partial M)} \quad\text{in }\Omega
\end{equation}
in the weak sense. Hence, by Theorem A.5 in \cite{FKS}, we have $|D^2 u|\in L_\loc^\infty(\Omega\setminus\Sigma)$. This, combined with Lemma 4.1 in \cite{FKS} and our assumption that $u\in\cN(\partial M)$ a.e.\ in $B$ yields, for each compactly contained open connected set $U\subset B\setminus\Sigma$, a constant $c > 0$ such that 
$$
|\Delta (\Pi\circ u)|\leq c |D(\Pi\circ u)| \quad\text{in }U. 
$$
Now, as we assume that $|D(\Pi\circ u)| =0$ on $F\subset B$ and that $|F|  >0$, the classical UCP due to Aronzjain (see Remark 3 in \cite{A}) applies, therefore $|D(\Pi\circ u)| = 0$ in the maximal open connected set in $B\setminus\Sigma$ containing $F$. However, since the set $B\setminus\Sigma$ itself is connected by assumption, and since $|D(\Pi\circ u)|\in L^2(B)$, we must have $|D(\Pi\circ u)| = 0$ in $B$. In other words, we have
\begin{equation}
    \label{eq:DPu-0}
    \Pi\circ u = y_0 \quad\text{in }B, 
\end{equation}
for some $y_0\in\partial M$. Noting that $u = \Pi\circ u$ in $u^{-1}(\partial M)$, it follows from \eqref{eq:u-pde} and \eqref{eq:DPu-0} that 
\begin{equation}
    \label{eq:u-pde-re}
    \Delta u = 0 \quad\text{in }B;
\end{equation}
in particular, $B\cap\Sigma = \emptyset$. By \eqref{eq:DPu-0} and \eqref{eq:u-pde-re}, we must have 
\begin{equation}
    \label{eq:u-re}
    u = y_0 + h\nu(y_0)\quad\text{in }B,
\end{equation}
for some harmonic function $h$ on $B$ which, because of the constraint $u(x) \in \overline M$ a.e., satisfies $h\geq 0$.

Now, by the strong minimum principle, we must have either $h = 0$ in $B$ or $h > 0$ in $B$. Let us treat these two cases separately.

\begin{claim}
    \label{claim:case1}
    If $h = 0$ in $B$, then $u = y_0$ in $\Omega$. 
\end{claim}

Indeed, if $h = 0$ in $B$, then \eqref{eq:u-re} yields $u = y_0$ in $B$. Since $|D^2u|\in L_\loc^\infty(\Omega\setminus\Sigma)$, \eqref{eq:u-pde} yields  
\begin{equation*}
    \label{eq:u-pde-re2}
    |\Delta u|\leq c|Du|\quad\text{in }U
\end{equation*}
for any open set $U \subset \overline U \Subset \Omega\setminus\Sigma$, with a constant $c$ depending  on $U$. Thus, $u$ satisfies the UCP in $U$, which implies that $u = y_0$  in the connected set $\Omega\setminus\Sigma$ containing $B$,
as desired.

\begin{claim}
    \label{claim:case1}
    If $h > 0$ in $B$, then it can be extended to a harmonic function in the entire domain $\Omega$, and $u = y_0 + h
    \nu(y_0)$ in $\Omega$.  
\end{claim}

Let $V\subset\Omega\setminus\Sigma$ be the maximal open set containing $B$ for which 
$$
u(V)\subset \{y_0 + t\nu(y_0):t>0\}.
$$
We shall prove that $\overline V\cap\Omega\setminus\Sigma\subset V$. Then $V$ being both open and (relatively) closed must be the entire $\Omega\setminus\Sigma$.

Suppose towards a contradiction that $\partial V\cap (\Omega\setminus\Sigma)\neq\emptyset$. Then there is a point $x_0\in\partial V\cap\Omega\setminus\Sigma$ and, since \eqref{eq:u-re} holds in $V$, by definition of $V$ we must have that $u(x_0) \in \partial M$. Moreover, since $x_0\not\in\Sigma$, we can find some small $\delta > 0$ for which $u(B_\delta(x_0))\subset \cN(\partial M)$. Since by assumption $\{y_0 + t\nu(y_0):t>0\}\cap M^c=\emptyset$, we must have $y_0 = u(x_0)$. Thus, $\Pi\circ u = y_0$ in $B_\delta(x_0)\cap V$. This puts us back to our initial setting, but now with $B = B_\delta(x_0)$ and $F = \overline{B_\delta(x_0)\cap V}$. Arguing exactly as in the derivation of \eqref{eq:u-re}, we obtain 
$$
u = y_0 + \tilde h \nu(y_0) \quad\text{in }B_\delta(x_0),
$$
for some harmonic function $\tilde h \geq 0$ in $B_\delta(x_0)$. However, as $u(x_0) = y_0$, we must have $\tilde h(x_0) = 0$, which then by the minimum principle implies $\tilde h = 0$ in $B_\delta(x_0)$. This is a contradiction to $B_r(x_0)\cap V \neq\emptyset$, proving that $V = \Omega\setminus\Sigma$.

Recalling that $\{y_0 + t\nu(y_0):t>0\}\cap M^c=\emptyset$, we must have $u(\Omega\setminus\Sigma) \subset M$, so it follows from \eqref{eq:u-pde} and $|\Sigma| = 0$ that 
$$
\Delta u = 0\quad\text{in }\Omega. 
$$
This in turn implies $\Sigma = \emptyset$, so $V= \Omega$ and the desired conclusion follows. 
\end{proof}

Building upon the preceding result, the following corollary shows  that if $\Pi\circ u$ is constant on the boundary of an interior ball, then it remains constant throughout, provided the obstacle satisfies condition \eqref{eq:geometriccond}.

\begin{cor}\label{cor:proj-ucp}
Let $u\in W^{1,2}(\Omega;\overline M)$ be a minimizing constraint map such that $u\in \cN(\partial M)$ a.e.\ in a ball $B\subset\Omega$, and that $\Pi\circ u = y_0$ on $\partial B$ for some point $y_0\in\partial M$. If $M$ satisfies \eqref{eq:geometriccond}, then $u=y_0+h\nu(y_0)$ in $\Omega$ for some non-negative harmonic function $h$ in $\Omega$, and either $h>0$ or $h\equiv 0$. 
\end{cor}

 \begin{proof}
 Define $\vp: \Omega\to\R^m$ by 
     $$
     \vp =\begin{cases}
         y_0 + (\rho\circ u)\nu(y_0)& \text{in }B,\\
         u & \text{in }\Omega\setminus B. 
     \end{cases}
     $$
 Since $M$ satisfies \eqref{eq:geometriccond} and $\Pi\circ u- y_0 \in W_0^{1,2}(B;\R^m)$, 
 we deduce that $\vp \in W^{1,2}(\Omega;\overline M)$ with  $\supp(\vp - u) \subset B\Subset \Omega$: in other words, $\vp$ is an admissible map. Thus the minimality of $u$ implies that 
 $$
 \int_B |Du|^2\,dx \leq \int_B |D\vp|^2\,dx = \int_B |D(\rho\circ u)|^2\,dx.
 $$
 Since $u= (\rho \circ u) \nu\circ u + \Pi \circ u,$ and as $|\nu|=1$ and $(\nabla\Pi)\nu=0 $, we have
 $$|Du|^2 = |D(\rho \circ u)|^2 + |(\rho \circ u) D(\nu\circ u)+D(\Pi\circ u)|^2 = |D(\rho \circ u)|^2 + |(Du)^\top|^2,$$
 where $(Du)^\top$ is the projection of $Du$ onto $T(\partial M)$, cf.\ \cite[(2.4)]{FGKS2}. 
 
 This implies that $|(Du)^\top| =0$ a.e.\ in $B$, and since $u\in \cN(\partial M)$ a.e.\ in $B$, we infer that also $D(\Pi \circ u)= \cA_u(Du)^\top=0$ a.e.\ in $B$ (here $\cA_y$ is an invertible matrix whenever $y\in\cN(\partial M)$, cf.\ \cite[(2.4) and (2.10)]{FGKS2}, and we wrote $\xi^\top = (\xi_\alpha^\top)_\alpha$ and $\cA_y\xi^\top = (\cA_y\xi_\alpha^\top)_\alpha$). Note also that the energy minimality of $u$ implies that it is continuous in $\Omega$ away from a relatively closed set $\Sigma$ of Hausdorff dimension at most $n-3$, according to \cite{DF}. In particular, $\Sigma$ cannot disconnect any open-connected set in $\Omega$. Then, the conclusion of the corollary follows from Proposition \ref{prop:proj-ucp}. 
\end{proof}

Corollary \ref{cor:proj-ucp} is known in the case of smooth harmonic maps (not necessarily minimizing) into closed manifolds, see e.g.\ \cite{Lemaire1978} for $n=2$ and \cite{KW} for general dimension. It may be that Corollary \ref{cor:proj-ucp} also holds for sufficiently regular weakly constraint maps, or for more general obstacles, but we do not pursue this here, leaving further exploration on this matter to future research, cf.\ Problem \ref{prob:constantPi}.

We note that, as asserted in Theorem \ref{thm:UCPPi}, the geometric condition \eqref{eq:geometriccond} on $M$ in Proposition~\ref{prop:proj-ucp} is optimal: the next example shows a situation where there is $a\in\partial M$ such that the ray $\{a + t\nu(a): t > 0\}$ touches $M^c$ tangentially at some other point $b\in\partial M$ and the conclusion of the proposition does not hold.

\begin{ex}\label{ex:proj-ucp}
Let $n=1, m=2$ and consider the obstacle $M^c$ depicted in Figure \ref{fig:Cobstacle}, which does not satisfy the geometric condition \eqref{eq:geometriccond}. Consider as map $u\colon [0,1]\to \overline M$ the unit-speed parametrization of the unique minimizing geodesic in $\overline M$ connecting $u(0)$ and $u(1)$, as depicted. Clearly there is $\delta>0$ such that $u([0,\delta])\subset \cN(\partial M)$ and furthermore $\Pi\circ u$ is constant in $[0,\delta]$. Yet $\Pi\circ u$ is not constant throughout $[0,1]$ and $u$ is not of the form $u=y_0 + h\nu(y_0)$ for some $y_0\in \partial M$ and some harmonic (i.e.\ affine) function $h\colon [0,1]\to [0,+\infty)$.
\begin{figure}[h]
    \centering
    \includegraphics[scale=0.75]{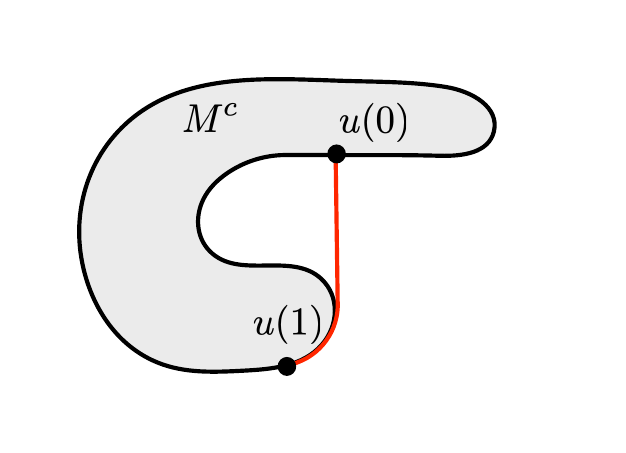}
    \caption{An obstacle which does not verify \eqref{eq:geometriccond}.}
    \label{fig:Cobstacle}
\end{figure}

\end{ex}

Proposition \ref{prop:proj-ucp} and Example \ref{ex:proj-ucp} suggest that, in general, rank one gradients of constraint maps may satisfy the UCP: if $u$ is a regular constraint map with $\rank Du = 1$ in a ball $B\subset\Omega$, then $\rank Du = 1$ in an open dense subset of $\Omega$, and $u(\Omega)$ lies on a geodesic arc in $\overline M$. This assertion is true for harmonic maps, see \cite[Theorem 3]{S} and \cite{HL}, but it (surprisingly) fails for higher rank unless the metric is analytic \cite{Jin1991}. Indeed,  there are smooth harmonic maps for which $\rank(Du)=2$ and $\rank(Du)=3$, respectively, in two disjoint non-empty open sets.
The next example shows that the UCP property of rank-one gradients fails, in general, for constraint maps.

\begin{ex}\label{ex:ucprank}
Let $n=m=2$, and consider an obstacle $M^c$ as in Figure \ref{fig:valley}; note that $M$ can be arranged to satisfy \eqref{eq:geometriccond}. We take $\Omega=\mb B^2$ to be the unit disk, we prescribe boundary conditions as depicted, 
and we let $u$ be a corresponding minimizing constraint map. 
   \begin{figure}[h]
    \centering
    \includegraphics[scale=0.4
    ]{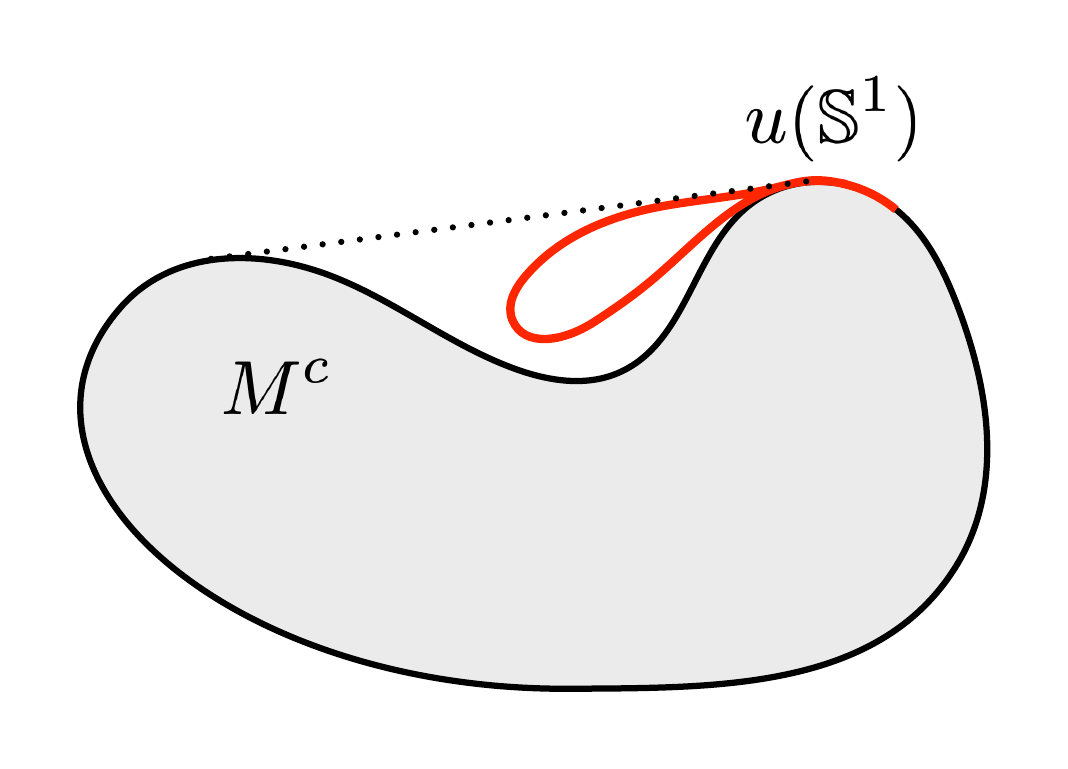}
    \caption{An obstacle and boundary data that force $\textup{rank}(Du)$ to be 1 in some places and 2 in others.}
    \label{fig:valley}
\end{figure}
More precisely, suppose that the boundary condition maps the arc segment $\{e^{i \theta}:\theta \in [-\frac \pi 2, \frac \pi 2]\}\subset \mb S^1$ into $\partial M$, and also that $\dist(u|_{\mb S^1}, Y)=0$, where $Y$ is the closed convex hull of $M^c$. Since $\dist(u,Y)$ is subharmonic (cf.\ \cite[Lemma 3.2]{FGKS2} or Lemma \ref{lem:dist-cont} below) we infer that $\dist(u,Y)=0$ in $\mb B^2$. Then, since $u\in C^1(\overline{\mb B^2})$ as $n=2$, there is $r_0\in (0,1)$ such that $u(\mb A)\subset \partial M$, where $\mb A=\{re^{i \theta}: r\in (r_0,1), \theta \in (-\frac \pi 4, \frac \pi 4)\}$ is a suitable annular sector.
In particular, in $\mb A$ we necessarily have $\textup{rank}(Du)\leq 1$. 

Now suppose that the boundary condition maps the arc segment $\{e^{i \theta}:\theta \in [\frac {3\pi}{4}, \frac{5\pi}{4}]\}\subset \mb S^1$ into a non-straight curve in $M$. Again by continuity of $u$, we can find a neighborhood $U$ in $\overline {\mb B^2}$ of this arc segment such that $u(U)\subset M$. Since $\Delta u =0$ in $\textup{int}(U)$ we have that $u$ is analytic in $U$ and therefore $\textup{rank}(Du)=k$ in a dense open subset $U'$ of $U$. Clearly $k\in \{0,1,2\}$; we claim that $k=2$. If not, then $u(U')$ is a straight-line \cite[Theorem 3]{S}, hence by continuity also $u(U)$ is a straight line, contradicting our choice of boundary condition. 

This discussion shows that, inside its domain, $Du$ has rank 1 in some regions and rank 2 in others.

\end{ex}

\section{Unique continuation and regularity around the free boundary}\label{sec:ucpfb}

This section builds on the analysis in~\cite{FGKS, FGKS2} concerning the interplay between free boundaries and discontinuity points, and serves as motivation for the results in Section~\ref{sec:Ainfty}.
In order to simplify our exposition, we shall focus on the case of the Dirichlet energy considered in \cite{FGKS2}, leaving out the Alt--Caffarelli energy considered in \cite{FGKS}; in fact, the former turns out to be substantially more difficult than the latter when it comes to the issue we are addressing here. 

The statement of \cite[Theorem 1.1]{FGKS2}  asserts, in lay terms, that if the obstacle $M^c$ is uniformly convex then there is a universal neighborhood of the free boundary $\partial u^{-1}(M)$ where $u$ is free of discontinuous points. Moreover, $u$ admits a sharp {\it a priori} estimate there.
The proof of this important results combines a set of new ideas, including a (weaker) version of the results in Section \ref{sec:Ainfty} in the present paper. 

Let us show how the UCP comes into play in a simple scenario. Since it will play an important part in our proof, for the reader's convenience we begin by proving a more general version of Theorem \ref{thm:classic}, valid for maps in the class of Definition \ref{def:F}; the argument is similar to a part of the proof of Proposition \ref{prop:proj-ucp} in the previous section.

\begin{thm}[Classical UCP]\label{thm:UCPF}
    Let $u\in \cF(\Omega;\overline M)$,  defined in Definition \ref{def:F}.
    If $|Du|=0$ in a set of positive measure in $\Omega$ then $|Du|=0$ a.e.\ in $\Omega$.
\end{thm}

\begin{proof}
    Thanks to property (i) in Definition \ref{def:F} we know that $u$ satisfies the Euler--Lagrange system \eqref{eq:EL}, hence 
    $$|\Delta u|\leq c|D u|^2$$ 
    weakly in $\Omega$, for some constant $c > 0$; we can choose $c$ as the uniform bound for all principal curvatures of $\partial M$, which is possible as $\partial M$ is compact and smooth. Therefore, in any open connected set $U\subset\Omega$ in which $u$ is continuous, we deduce from the regularity theory for elliptic systems that $|D u|\leq C_U$ in $U$. Therefore we have
    $$|\Delta u|\leq cC_U|D u|$$ in $U$, so the UCP of \cite{A} applies to $u$. 

    Now, since by assumption
    $|Du|=0$ on a set $B\subset \Omega$ of positive measure, we deduce that $|D u| = 0$ in the maximal open connected set in $\Omega\setminus\Sigma$ containing $B$, where $\Sigma$ is the set where $u$ is discontinuous (as in (iii) of Definition \ref{def:F}). However, (iii) together with $|Du|\in L^2(\Omega)$ imply that $\dim_\cH\Sigma \leq n-2$, whence $\Omega\setminus\Sigma$ has to be connected. Thus, $|Du| = 0$ in $\Omega\setminus\Sigma$. Since $|Du|\in L_\loc^2(\R^n)$ and $|\Sigma| = 0$, the claim follows.
\end{proof}

We now continue with the following lemma, which is very close to previous results contained in \cite{FGKS,FGKS2}.

\begin{lem}
    \label{lem:dist-cont}
    Let $M$ be a smooth domain in $\R^m$ with convex and compact complement, and let $u\in \cF(\Omega;\overline M)$ be given. Then $\dist(u,M^c)$ is continuous and weakly subharmonic in $\Omega$. 
\end{lem}

\begin{proof}
    The composition of the convex function $\dist(\cdot,M^c)$ with a harmonic vector-valued function is subharmonic, hence $\dist(u,M^c)$ is subharmonic in $\textup{int}(\{\{\dist(u,M^c)>0\}).$ On the other hand, since $u$ is continuous outside a closed singular set $\Sigma$ with $\cH^{n-2}(\Sigma) =0$, by a capacity argument we conclude that $\dist(u,M^c)$ is subharmonic throughout $\Omega$.

    Since $\dist(u,M^c)$ is continuous in $\Omega\setminus\Sigma$, it suffices to check that $\dist(u,M^c)$ vanishes continuously on $\Sigma$. With properties (i)--(iv) from Definition~\ref{def:F} at our disposal, the proof for the latter then follows essentially the same line of arguments as that of Theorem 6.1 in \cite{FGKS}.
    To keep the discussion self-contained, let us show the argument here. 
    
    To see that $\dist(u,M^c)$ vanishes continuously on $\Sigma$, let us fix a point $x_0\in\Sigma$. In view of (ii) and (iii), there exists a constant $\e_0>0$ such that 
    $$
    e := \lim_{r\downarrow 0} r^{2-n}\int_{B_r(x_0)} |Du|^2\,dx \geq \e_0. 
    $$
    By subharmonicity it follows that $\dist(u,M^c)\in L_\loc^\infty(\Omega)$, so by the compactness of $M^c$ we deduce that $|u|\in L_\loc^\infty(\Omega)$. Thus, by (ii) and (iv), there exists a tangent map $\vp$, satisfying
\begin{equation}
    \label{eq:auxphi}
    R^{2-n}\int_{B_R} |D\vp|^2\,dx = e
\end{equation}
    for every $R > 0$. In particular, $\vp$ is $0$-homogeneous.

    Let $r_k\downarrow 0$ be the sequence for which $u_{x_0,r_k}\to \vp$ in $W^{1,2}(B_R;\R^m)$, for any $0< R < 1/r_k$. As $\vp\in \cF(B_R;\overline M)$ for every $R>0$ by (iv), the above reasoning yields $\Delta\dist(\vp,M^c) \geq 0$ in $\R^n$. Then the $0$-homogeneity of $\vp$ yields that $\dist(\vp,M^c) = a$ in $\R^n$, for some constant $a\geq 0$, as shown e.g.\ in \cite[Lemma 2.6]{FGKS}.

    If $a > 0$, then (i) applied to $\vp$ implies that $\Delta \vp = 0$ in $\R^n$. However, $\vp$ being a tangent map of a locally bounded map, it is bounded in the entire space $\R^n$: 
    $$\|\varphi\|_{L^\infty(\R^n)} = \sup_{r>0}\|\varphi\|_{L^\infty(B_{1/r})} = \sup_r \lim_{k\to \infty} \|u_{r_k}\|_{L^\infty(B_{1/r})}
    = \sup_r \lim_{k\to \infty}\|u\|_{L^\infty(B_{r_k/r})} \leq \|u\|_{L^\infty(B_1)}.$$
    Therefore, the Liouville theorem implies $|D\vp| = 0$ in $\R^n$. This however contradicts \eqref{eq:auxphi}.

    Thus, $\dist(\vp,M^c) = 0$ in $\R^n$. Then as $u_{x_0,r_k}\to \vp$ in $L_\loc^2(\R^n;\R^m)$, we have $\dist(u_{x_0,r_k},M^c)\to 0$ in $L_\loc^2(\R^n)$, whence for any $\e > 0$, we can find $k_\e\in\N$ large such that $\|\dist(u_{x_0,r_k},M^c)\|_{L^2(B_2)} \leq \e$ for all $k\geq k_\e$. Setting $\delta = r_{k_\e}$, we deduce from the local $L^\infty$-estimate for weakly subharmonic functions that 
    $$
    \| \dist(u,M^c) \|_{L^\infty(B_\delta(x_0))} \leq cr_k^{-\frac{n}{2}}\|\dist(u,M^c)\|_{L^2(B_{2\delta}(x_0))} \leq c\e. 
    $$ 
    As $\e > 0$ was arbitrary and $c > 0$ is independent of $\e$, we conclude that $\dist(u,M^c)$ vanishes continuously at $x_0\in\Sigma$, as desired. 
\end{proof}

Given $\delta>0$ we say that a set $E$ is $\delta$-porous in $U$ if, for every ball $B_r(x_0)\subset U$, we can find a smaller ball $B_{\delta r}(x)\subset B_r(x_0)\setminus E$. The constant $\delta$ is called the porosity constant. 

\begin{thm}
    \label{thm:porous}
    Let $M$ be a smooth domain in $\R^m$ with uniformly convex complement. Let $u\in \cF(\Omega;\overline M)$ be given, and assume that $\Omega\cap \partial u^{-1}(M)$ is $\delta$-porous in the closure of $u^{-1}(M)$. Then $u$ is continuous at any point of $\Omega\cap \partial u^{-1}(M)$.
\end{thm}

\begin{proof}
    Let $x_0\in\Omega\cap\partial u^{-1}(M)$ be given. Consider the tangent map $\vp\in \cF(\R^n;\overline M)$ of $u$ at $x_0$; the existence of a tangent map follows from property (ii) in Definition \ref{def:F} and the subharmonicity of $\dist(u,M^c)$. 
    Since $\dist(u(x_0),M^c) = 0$, we must have $\dist(\vp,M^c) = 0$ in $\R^n$, i.e., $\vp(\R^n)\subset\partial M$. Then by properties (iv) and (i) in Definition \ref{def:F}, we also have 
    \begin{equation}
        \label{eq:vp-pde}
        \Delta \vp = A_\vp(D\vp,D\vp)\quad\text{in }\R^n.
    \end{equation}
    Let $r_k\downarrow 0$ be a sequence along which $u_{x_0,r_k}\to \vp$ in $W_\loc^{1,2}(\R^n;\R^m)$ (the strong convergence follows from (iv)). By the $\delta$-porosity of $\Omega\cap\partial u^{-1}(M)$ in the closure of $u^{-1}(M)$, we find a sequence $x_k\in u^{-1}(M)$  such that $x_k\in B_{r_k}(x_0)$ and $B_{\delta r_k}(x_k)\subset B_{r_k}(x_0)\cap u^{-1}(M)$ for all $k\in\N$. In particular (i) implies that $\Delta u = 0$ in $B_{\delta r_k}(x_k)$.

Write $\xi_k := r_k^{-1}(x_k-x_0) \in B_1$, and assume without loss of generality that $\xi_{k_i}\to \xi$ for some $\xi\in B_1$ (with $B_\delta(\xi)\subset B_1$). Then we deduce from the $L^2$-convergence of $Du_{x_0,r_k}$ to $D\vp$ that $\Delta\vp = 0$ in $B_\delta(\xi)$. Comparing this with \eqref{eq:vp-pde}, we obtain $A_\vp(D\vp,D\vp) = 0$ in $B_\delta(\xi)$. As $M^c$ is uniformly convex, there is $\kappa > 0$ such that $|A_y(\xi,\xi)|\geq \kappa|\xi|^2$ for any $\xi\in T_y(\partial M)$ uniformly for all $y\in\partial M$. Thus, we arrive at 
    \begin{equation*}
        \label{eq:Avp}
        |D\vp| = 0\quad\text{in }B_\delta(\xi).
    \end{equation*}
    Since $\varphi \in \cF(\R^n;\overline M)$,  Theorem \ref{thm:UCPF} implies that 
    \begin{equation}
        \label{eq:Avp-re}
        |D\vp| = 0 \quad\text{in }\R^n.
    \end{equation}
    By \eqref{eq:Avp-re}, we obtain from property (iv) that $|Du_{x_0,r_k}| \to 0$ in $L_\loc^2(\R^n)$ as $r_k\downarrow 0$. Recalling that $\vp$ was an arbitrary tangent map of $u$ at $x_0$, we deduce that $|Du_{x_0,r}|\to 0$ in $L_\loc^2(\R^n)$ as $r\downarrow 0$, which in the non-rescaled form yields 
    $$
    \lim_{r\downarrow 0} r^{2-n}\int_{B_r(x_0)} |Du|^2\,dx = 0. 
    $$
    Recalling property (iii) in Definition \ref{def:F}, we conclude that $u$ is continuous at $x_0$, as desired. 
\end{proof}

With a slightly more careful argument (involving the PDE for the distance map $\dist(u,M^c)$), we can also extend Theorem \ref{thm:porous} to free boundary points having positive  density with respect to the non-coincidence set. We should also remark that all the above analysis goes through for constraint maps with respect to the Alt-Caffarelli energy \cite{FGKS}. For the latter work, the distance map a satisfies nondegeneracy condition, which ensures a uniformly positive density for all free boundary points. Consequently, we obtain full regularity of energy minimizing constraint maps near the non-coincidence set. 

However, the situation changes drastically in the case of the Dirichlet energy we consider here, as there is no inherent non-degeneracy of the distance map, see for instance Lemma 9.4 in \cite{FGKS2}. Thus, in this case it is unnatural to impose either positive density or a porosity condition of the free boundary with respect to the non-coincidence set. In \cite{FGKS2} we overcame the difficulty by establishing a quantitative version of the UCP, namely Theorem 1.2 therein. In the next section, we will prove a stronger version of that result.

\section{Unique continuation property and $A_\infty$-weights}\label{sec:Ainfty}

We begin by recalling the following standard definition, which the reader can find e.g.\ in \cite{Grafakos2009}.

\begin{defn}\label{def:Ainfty}
    Let $0\leq \omega\in L^1(B_4)$. We say that $\omega\in A_\infty(B_1)$ if there are constants $C,\e>0$ such that, for any ball $B$ with $2B\subset B_2$, the following hold:
    \begin{enumerate}
        \item weak reverse H\"older inequality:  
        $$\bigg(\fint_{B} \omega^{1+\e}(x) \, dx \bigg)^{\frac{1}{1+\e}}
        \leq C \fint_{2B} \omega(x) \, dx$$
        \item doubling: $\int_{2B} \omega(x) \, dx \leq C \int_{B} \omega(x) \, dx$.
    \end{enumerate}
    As usual, $2B$ denotes the ball concentric with $B$ but with twice the radius. We call $C,\e$ the $A_\infty$-constants of $\omega$.
\end{defn}

It is of course clear from the above definition that weights in $A_\infty$ are exactly those weights which satisfy a (standard) reverse H\"older inequality, i.e.\ a reverse H\"older inequality where the balls on the left- and right-hand sides can be chosen to be the same. One reason for defining $A_\infty$-weights as above is that, in applications to elliptic PDE, condition (1) is often a simple consequence of a Caccioppoli estimate (see already Corollary \ref{cor:revHolder} below), while it is typically the doubling condition (2) that is much more non-trivial to check.

There are several other equivalent characterizations of $A_\infty$-weights, and here we chose the one which is more natural to verify in applications to PDEs. The following theorem gives another characterization of $A_\infty$-weights:

\begin{thm}\label{thm:ucpA_infty}
    Let $0\leq \omega\in L^1(B_4)$. Then
    $\omega\in A_\infty(B_1)$ if and only if either $\omega=0$ a.e.\ or else
    there are $\tilde C,\gamma>0$ such that, whenever $2B\subset B_2$, we have
    \begin{equation}
        \label{eq:Ainftyequiv}
        \bigg(\fint_B \omega(x) \, dx \bigg)\bigg(\fint_{B} \omega^{-\gamma}(x) \,dx \bigg)^{-\frac 1 \gamma} \leq \tilde C.
    \end{equation}
    This result is quantitative in the sense that the constants $\tilde C, \gamma$ depend only on $C,\e$ from Definition \ref{def:Ainfty} and on $n$. Conversely, if \eqref{eq:Ainftyequiv} holds, then the $A_\infty$-constants of $\omega$ depend only on $\tilde C, \gamma$ and $n$. 
\end{thm}

Theorem \ref{thm:ucpA_infty} shows that $A_\infty$-weights have a very strong, quantitative form of the UCP: either they vanish or else they detach from 0 (in an integral sense) at a universal rate. Theorem \ref{thm:ucpA_infty} is well-known, see e.g.\ \cite{Coifman1974,Grafakos2009}. To be precise, in these references the result is stated in a global manner (i.e.\ for weights defined over $\R^n$), but it is easy to see that the proofs are local in nature, cf.\ the appendix in \cite{DePhilippis2023}.

Using Theorem \ref{thm:ucpA_infty}, we see that Theorem 1.2 in \cite{FGKS2} is in particular implied by the following:

\begin{thm}\label{thm:L-gamma}
    Let $N\subset \R^m$ be a smooth, compact manifold without boundary and $u\in W^{1,2}( B_4; N)$ be a minimizing harmonic map. Then $|Du|\in A_\infty(B_1)$. In fact, if $u$ is non-constant then there is $y_0\in B_2$ and constants $0<\lambda<\Lambda$ such that
	\begin{equation}
	\label{eq:assumu}
	0<\lambda \leq \int_{B_1(y_0)}|Du|^2 \, dx\leq \int_{B_4}|Du|^2 \, dx\leq \Lambda,
	\end{equation}
    and the $A_\infty$-constants of $|Du|$ depend only on $\Lambda,\lambda,n,N$.
\end{thm}

The conclusion of Theorem \ref{thm:L-gamma} holds more generally for suitable manifolds $N$ with boundary. In Remark \ref{rem:manifoldswithbdry2} we will indicate how to modify the proof in this case.

To  prove  Theorem \ref{thm:L-gamma} we need a few technical lemmas.
We will use without further comment that minimizing harmonic maps satisfy properties (i)-(iv) in Definition \ref{def:F}; in fact, we also recall that (iii) can be upgraded to the statement that 
\begin{equation}
\label{eq:singset}
\Sigma(u) = \left\{x\in \Omega:\lim_{r\to 0} E(u,x,r) \geq \e_0^2\right\}
\end{equation}
for some $\e_0=\e_0(n,N)>0$.
 We will need an additional fundamental property of minimizing harmonic maps, which is closely related to the compactness property (iv): they satisfy a Caccioppoli inequality, see e.g.\ \cite[Lemma 1, \S 2.8]{Simon}.

\begin{lem}\label{lem:caccioppoli}
    Let $u\in W^{1,2}( B_2; N)$ be a minimizing harmonic map with $\int_{B_2} |Du|^2 \, dx \leq \Lambda$. There is a constant $C$, depending only on $\Lambda,n,N$, such that, for any ball $B_r(x_0)\subset B_1,$
    $$\int_{B_{r/2}(x_0)} |Du|^2 dx \leq C r^{-2}\int_{B_r(x_0)} |u-u_{x_0,r}|^2 dx,$$
    where $u_{x_0,r}=\fint_{B_r(x_0)} u\, dx$.
\end{lem}
    
As an immediate consequence, we obtain the following corollary, which is well known to the experts (see e.g.\ \cite{Hardt1988,Hardt1987a} for similar results when the target $N$ is simply connected):

\begin{cor}\label{cor:revHolder}
    Let $u\in W^{1,2}( B_2; N)$ be a minimizing harmonic map with $\int_{B_2} |Du|^2 \, dx \leq \Lambda$. There is $p>2$ and $C>0$, depending only on $\Lambda,n,N$, such that, for any ball $B_r(x_0)\subset B_1,$
    $$\bigg(\fint_{B_{r/2}(x_0)} |D u|^p \, dx \bigg)^{\frac 1 p }\leq C \bigg(\fint_{B_r(x_0)} |D u|^2 \, dx \bigg)^{\frac 1 2 }.$$
\end{cor}

\begin{proof}
 Applying Lemma \ref{lem:caccioppoli} and the Poincar\'e--Sobolev inequality, we have
$$\bigg(\int_{B_{r/2}(x_0)} |Du|^2\, dx\bigg)^{\frac 1 2} 
\leq C \bigg(\int_{B_{r}(x_0)} |Du|^{2_*} \, dx\bigg)^{1/2_*},$$
where $2_*=\frac{2n}{n+2}<2$. Hence $Du$ satisfies a weak reverse H\"older inequality. By the general theory of such reverse H\"older inequalities, see e.g.\ \cite[\S 14]{Iwaniec2001}, the conclusion follows.
\end{proof}

Corollary \ref{cor:revHolder} shows that, to prove Theorem \ref{thm:L-gamma}, we only need to check the doubling condition.
We note in passing that 
the higher integrability result in Corollary \ref{cor:revHolder} 
yields the compactness of minimizing harmonic maps, since it forbids energy concentration, see e.g.\ \cite[\S 4]{LW}. We also remark that it is known from the works \cite{CN,NV1} that in fact $Du\in \textup{weak-}L^3$, although we will not need this much deeper fact here. 

As a consequence of higher integrability, we obtain:

\begin{cor}\label{cor:noconcen}
Let $u\in W^{1,2}(B_2;N)$ be a minimizing harmonic map with $\int_{B_2} |Du|^2 \, dx \leq \Lambda$. There are constants $\theta,C>0$, depending only on $\Lambda,n,N$, such that
    $$\int_{B_{1+\delta}\setminus B_1} |Du|^2 \, dx \leq C\delta^\theta \int_{B_2} |D u|^2 \, dx$$
\end{cor}

\begin{proof}
    Given $B_{2r}(x_0) \subset B_2$, by Corollary \ref{cor:revHolder} and H\"older's inequality we have, with $\theta= \frac{n(p-2)}{p}>0,$
    $$\int_{B_{r}(x_0)} |Du|^2 \, dx\leq C r^{\theta}\bigg(\int_{B_{2r}(x_0)} |Du|^p \, dx\bigg)^{\frac 2 p}\leq C r^\theta \int_{B_{2r}(x_0)} |Du|^2 \, dx.$$
    Covering the annulus $B_{1+\delta}\setminus B_1$ by $C(n)\delta^{1-n}$ balls of radius $\delta\leq \frac 1 4$, we have
    $$\int_{B_{1+\delta}\setminus B_1} |Du|^2 \, dx \leq C \delta^{\theta}\int_{B_{1+4\delta}\setminus B_{1-2\delta}} |Du|^2 \, dx \leq C \delta^\theta \int_{B_2} |Du|^2 \, dx$$
    as wished.
\end{proof}

Continuing with the proof of Theorem \ref{thm:L-gamma}, we introduce the following definition from \cite{FGKS2}:

\begin{defn}
Let $u \in W^{1,2}(B_4,N)$ be a minimizing harmonic map. Given $\ell_0\geq 2$, we set
$$r_u(x) := \begin{cases}
\sup\left\{ r\in (0,2]: E(u,x,r)\leq \frac{\e_0^2}{\ell_0^2}\right\} & \text{if }x\in B_2\setminus \Sigma(u),\\
0 &  \text{if }x\in B_2\cap \Sigma(u),
\end{cases}$$
to be the \textit{critical scale} of $u$ at $x$.
\end{defn}

The definition of critical scale depends on $\ell_0\geq 1$, but a simple compactness argument shows that in fact this parameter can be chosen uniformly:

\begin{lem}\label{lem:critscale}
Let $u$ be as in Theorem \ref{thm:L-gamma} and  suppose that \eqref{eq:assumu} holds.  There is $\ell_0=\ell_0(\lambda,\Lambda,n,N)\geq 2$ such that $\Sigma(u) = r_u^{-1}(0)$ and $r_u\leq 1$ in $B_2$.
\end{lem}

\begin{proof}
This is essentially shown in \cite[Lemma 4.1]{FGKS2}, so we just sketch the main points. 

First, recalling \eqref{eq:singset}, it follows that $\Sigma(u) = r_u^{-1}(0)$. 

The bound $r_u\leq 1$ is proved by contradiction: if it were to fail, using the compactness of the family of minimizing harmonic maps with bounded energy, we would find a minimizing harmonic map $u\in W^{1,2}(B_3,N)$ and points $x_0,y_0\in \overline B_2$ such that 
$$E(u,x_0,1)=0,\qquad E(u,y_0,1)\geq \lambda>0.$$
This is a contradiction to the UCP of minimizing harmonic maps.
\end{proof}

In the rest of this section, we choose $\ell_0$ so that the conclusion in Lemma \ref{lem:critscale} holds.

The critical scale separates scales where $u$ behaves like a solution to a regular elliptic system from scales where the energy of $u$ concentrates. One of the main observations in \cite[\S 4]{FGKS2} is that, at the critical scale, the frequency is universally bounded, where we recall that the frequency is defined as
$$\mathscr{N}  (u,x_0,r):=\frac{r \int_{B_r(x_0)} |Du|^2\, dx }{\int_{\partial B_r(x_0)} |u-(u)_{x_0,r}|^2 \, d \cH^{n-1}}.$$
 Here we prove the following variant of 
\cite[Lemma 4.2]{FGKS2}:

\begin{lem}[Frequency bound above the critical scale]  \label{lem:freqbound}
Let $u$ be as in Theorem \ref{thm:L-gamma} and suppose that \eqref{eq:assumu} holds. There is $\mathscr{N}_0=\mathscr{N}_0(\lambda,\Lambda,n,N)\geq 1$ such that 
    $$\mathscr{N}(u,x,r)\leq \mathscr{N}_0, \quad \text{whenever } x \in B_1 \text{ and } r_u(x)< r\leq 1.
    $$
\end{lem}

\begin{proof}
If the conclusion were false, we would find a sequence of minimizing harmonic maps $(u_k)\subset W^{1,2}(B_4,N)$, satisfying \eqref{eq:assumu} at suitable points $(y_k)\subset B_2$, but for some $(x_k)\subset B_1\setminus \Sigma(u_k)$ we would have
\begin{equation}
\label{eq:contrafreq}
\mathscr{N}(u_k,x_k,r_k)\geq k \quad \text{for } r_{u_k}(x_k)< r_k\leq 1.
\end{equation}
Writing $v_k(y):=u_k(\tfrac 1 2 r_k y + x_k)$, we have that $(v_k)\subset W^{1,2}(B_3, N)$ is a sequence of minimizing harmonic maps that, by the monotonicity formula, has bounded energy. Thus, passing to subsequences, there is another minimizing harmonic map $v\in W^{1,2}(B_2,N)$ such that $v_k \to v$ in $W^{1,2}(B_2,\R^m)$, and hence in particular also $v_k\to v$ strongly in $L^2(\partial B_2, \R^m)$.  Since $\mathscr{N}(v_k,0,2)=\mathscr{N}(u_k,x_k,r_k)\geq k$ and $v_k$ have bounded energy, the numerator in the definition of $\mathscr{N}(v_k,0,2)$ stays bounded, therefore the denominator must go to zero. We infer that 

$$0 = \lim_{k\to \infty} \int_{\partial B_2} |v_k-(v_k)_{0,2}|^2 \,d\cH^{n-1} = \int_{\partial B_2} |v-(v)_{0,2}|^2 \,d\cH^{n-1} ,$$
i.e.\ $v$ is constant on $\partial B_2$. The minimality of $v$ then  yields that $v$ must be constant in $B_2$ (cf.\ the proof of Corollary \ref{cor:proj-ucp}), and then by the UCP  in fact $v$ is constant in $B_3$. This gives us the desired contradiction, since
$$E(v,0,2) = \lim_{k \to \infty} E(v_k,0,2)= \lim_{k\to \infty} E(u_k,x_k,r_k) \geq \frac{\e_0^2}{\ell_0^2}$$
by the definition of the critical scale.
\end{proof}

\begin{proof}[Proof of Theorem \ref{thm:L-gamma}]
    We can suppose that $u$ is not constant. In particular, \eqref{eq:assumu} holds for some $y_0\in B_2$ and some suitable constants $0<\lambda<\Lambda$. 
    
    By Corollary \ref{cor:revHolder}, it suffices to prove the doubling property of $|Du|$, namely
    \begin{equation}
        \label{eq:doubling}
        \int_{B_{2r}(x_0)}|Du|^2\, dx \leq C \int_{B_r(x_0)}|Du|^2\, dx
    \end{equation}
    whenever  $x_0\in B_1$ and $r\leq 1$.

We first claim that \eqref{eq:doubling} holds if $r_u(x_0)\leq r \leq 1$. If not, then  there are sequences of minimizing harmonic maps $(u_j)\subset W^{1,2}(B_4;N)$, satisfying \eqref{eq:assumu} at suitable points $(y_k)\subset B_2$, and sequences $(x_j)\subset B_1$ and $r_{u_j}(x_j)\leq r_j\leq 1$ such that
    $$ j \int_{B_{r_j}(x_j)}|Du_j|^2 \, dx  \leq \int_{B_{2r_j}(x_j)} |D u_j|^2\, dx .$$
By definition of the critical scale and the monotonicity formula, this would yield
$$j \frac{\e_0^2}{\ell_0^2} \leq j  E(u_j,x_j,r_j) \leq 2^{n-2} E(u_j,x_j,2r_j) \leq 2^{n-2} E(u_j,x_j,2)\leq \int_{B_3} |Du|^2\, dx \leq \Lambda,$$
impossible.

We now prove \eqref{eq:doubling} when $0<r\leq \frac 1 2 r_u(x_0)$.  Let $v_0:=u(x_0+r_u(x_0)\cdot)$ and note that, by definition of the critical scale,
$$E(v_0,0,1) = E(u,x_0,r_{u}(x_0)) = \frac{\e_0^2}{\ell_0^2}< \e_0^2.$$
Clearly $v_0\in W^{1,2}(B_3;N)$ is a minimizing harmonic map with energy at most $\Lambda$ and hence, applying Corollary \ref{cor:noconcen}, we can find
$\theta_0=\theta_0(\e_0, \ell_0, \Lambda,n,N)\in (0,1)$ such that
$$E(v_0,0,1+\theta_0)\leq \e_0^2.$$
Thus, by the $\e$-regularity theorem, $v_0$ is a solution of a regular elliptic system of the form
$$\Delta v_0^i = a_{ij}^\alpha \partial_\alpha v_0^j, \qquad \|a_{ij}^\alpha\|_{L^\infty(B_{1+\theta_0/2})}\leq C.$$
Also, by Lemma \ref{lem:freqbound}, we have $\mathscr{N}(v_0,0,1+\frac{\theta_0}{2} )=\mathscr{N}(u,x_0,(1+\frac{\theta_0}{2})r_u(x_0))\leq \mathscr N_0$.
Here both $C$ and $\mathscr N_0$ depend only on our fixed parameters $\lambda,\Lambda,n,N$.
Thus, by frequency-based theory for such systems (see e.g.\ \cite[Appendix B]{FGKS2} and \cite{Garofalo1986,GL,HL}), there is a new constant $C$, with the same dependence, such that
$$\int_{B_{2s}}|Dv_0|^2\, dx  \leq C\int_{B_s}|Dv_0|^2\, dx  $$
for any $s\leq \frac 1 2$. Thus, rescaling back to $u$, we find that
$$\int_{B_{2r}(x_0)} |D u|^2 \, dx \leq C \int_{B_r(x_0)} |Du|^2 \, dx,$$
whenever $0<r\leq \frac 1 2 r_u(x_0)$. 

Finally, the fact that \eqref{eq:doubling} holds for $\frac 1 2 r_u(x_0)\leq r \leq r_u(x_0)$ clearly follows from what we have shown so far, since
\begin{align*}
    \int_{B_{2r}(x_0)} |Du|^2 \, dx 
     \leq \int_{B_{2r_u(x_0)}(x_0)}|Du|^2 \, dx & \leq C\int_{B_{r_u(x_0)}(x_0)}|Du|^2 \, dx \\& \leq C^2 \int_{B_{r_u(x_0)/2}(x_0)}|Du|^2 \, dx\leq C^2 \int_{B_{r}(x_0)}|Du|^2 \, dx.
\end{align*}
This completes the proof.

\end{proof}

\begin{rem}\label{rem:manifoldswithbdry2}
    The conclusion of Theorem \ref{thm:L-gamma} holds if $N=\overline M$ is the complement of a smooth, bounded domain in $\R^m$; in this case, however, one needs to replace the Dirichlet energy in \eqref{eq:assumu} with the full $W^{1,2}$-norm, since the maps do not have a universal $L^\infty$-bound anymore. More generally, Theorem \ref{thm:L-gamma} holds provided that $N$ is a complete manifold with boundary satisfying the conditions in Remark \ref{rem:manifoldswithbdry}. The proof in the general case is exactly the same as the one we presented here: the reader will note that, apart from the properties in Definition \ref{def:F}, we only used two other properties of minimizing harmonic maps, namely that $|\Delta u|\leq c|Du|$ away from the singular set $\Sigma(u)$, and that $u$ satisfies the Caccioppoli inequality in Lemma \ref{lem:caccioppoli}. The first property is true for general smooth targets $N$, cf.\ \cite{D}, while the Caccioppoli inequality follows by modifying the arguments used in the boundaryless case, using Luckhaus' Lemma \cite{L}.
\end{rem}

\section{Future directions and open problems}
\label{sec:probs}

In this section, we discuss possible developments and new directions in the emerging field of vectorial free boundary problems.

The theory of free boundary problems has branched into very diverse directions and different properties of solutions have been extensively studied. Many of these avenues and problems hold promise for applications to constraint maps. However, beyond the classical free boundary theory, there are also several other questions and directions that naturally emerge and are more unique to the study of constraint maps.

Two key challenges which are absent in most (if not all) classical free boundary problems are the possible lack of continuity of solutions (maps in the vectorial setting) and the nonlinear influence of the constraint, which appears as source term in \eqref{eq:EL} and depends on the derivative of the map, whose smoothness is a priori unknown. These aspects have recently been examined in the  papers \cite{FGKS, FGKS2}. The first difficult point is  establishing the continuity of solutions, which serves as the cornerstone for further developments of the theory. This point may prove significantly more challenging in certain variations of the  problem. For instance, the approach outlined in this note, based on the (quantitative) UCP, would be very difficult to implement if we study maps minimizing the $p$-Dirichlet energy for $p\neq 2$: recall that, even in the scalar case, the qualitative UCP is not known for $p$-harmonic functions, unless the domain is two-dimensional, where quasiconformal techniques become applicable.

A further, broad direction of research is to move beyond the study of minimizers.
In fact, in the scalar setting, several free boundary problems emerge as non-variational models\footnote{It is noteworthy that some problems may have both variational and non-variational formulations. In this case, the solutions are in general non-unique. The prototype example is the Bernoulli free boundary problem.}.
These include the Bernoulli free boundary problem, the non-variational Stefan problem, the incompressible Euler equations with free boundaries, the Muskat problem (two-phase fluid flow in porous media), the plasma problem (Prandtl-Batchelor free boundary problem), the no-sign obstacle problem, and the free boundary problem in the superconductivity model, among others. 

For systems and mappings, several model problems arise as non-variational problems, introducing new types of challenges in their analysis. Moreover, they often lead to entirely new questions and phenomena that do not typically appear in scalar settings. A concrete example, which was briefly considered in \cite{FKS}, is to study \textit{leaky maps}, which correspond to the vectorial version of solutions to the no-sign obstacle problem.

\medskip

For the readers interested in this fascinating research topic, we conclude this article by presenting a concise list of concrete open problems that we find particularly interesting. 

We begin our list by stating six problems of a general nature, i.e.\ they pertain to regularity and unique continuation properties of constraint maps, in the setting described in this paper.

\begin{prob}[Optimal conditions for $C^{1,1}$-regularity]\label{prob:optimalM}
What is the optimal regularity on $M$ which ensures that a minimizing constraint map $u\in W^{1,2}(\Omega;\overline M)$ is in $C^{1,1}_\loc(\Omega\setminus \Sigma(u))$?
\end{prob}

In \cite{FGKS2,FKS} we proved $C^{1,1}$-regularity of $u$ around regular points under the assumption that $\partial M\in C^{3,\alpha}$ for some $\alpha>0$, while $W^{2,p}$-regularity for $p<\infty$ only requires that $\partial M \in C^2$. Thus Problem \ref{prob:optimalM} asks whether this gap can be closed. A reasonable guess is that perhaps the optimal regularity of $u$ holds even for obstacles of class $\partial M\in C^{1,1}$.

\begin{prob}[Optimal regularity of the projection map]\label{prob:Pi}
Let $u\in W^{1,2}(\Omega;\overline M)$ be a minimizing constraint map. Is it the case that $D^3(\Pi\circ u)\in L^\infty_\loc(\Omega\setminus \Sigma(u))$? \end{prob}

We note that Problem \ref{prob:Pi} asks for the best possible regularity for $\Pi\circ u$, since in general $\Pi \circ u\not \in C^3(\Omega\setminus \Sigma(u))$. We note that in \cite{FKS} it was shown that $D^3(\Pi\circ u)\in L^\infty_\loc(\Omega\setminus \Sigma(u))$  when $n=2$, while for $n\geq 3$ we only have $D^3(\Pi\circ u)\in BMO_\loc(\Omega\setminus \Sigma(u))$.

\medskip

\begin{prob}[Constant $\Pi$ on the boundary]\label{prob:constantPi}
    Let $u\in W^{1,2}(B;\overline M)$ be a $C^{1,1}$ weakly constraint map. If $\Pi \circ u$ is constant on $\partial B$, is it the case that $\Pi\circ u$ is constant in $B$, at least when $M^c$ is convex? What are the optimal geometric conditions on $M$ for this to hold?
\end{prob}

We note that Corollary \ref{cor:proj-ucp} provides a partial answer to Problem \ref{prob:constantPi}. In the case where $\rho\circ u$ vanishes identically, the answer is affirmative \cite{KW,Lemaire1978}.

\medskip

\begin{prob}[UCP and low rank]\label{prob:UCPrank}
    Let $M^c$ be convex and $u\in W^{1,2}(\Omega;\overline M)$ be a minimizing constraint map. If $\textup{rank}(Du)\leq 1$ in an open subset of $\Omega$, do we have $\textup{rank}(Du)\leq 1$ in $\Omega$?
\end{prob}

Problem \ref{prob:UCPrank} is closely related to the UCP of the rank of a harmonic map \cite{S}, but Example \ref{ex:ucprank} shows that the situation is more subtle for constraint maps. In particular, additional assumptions on $M^c$ (such as convexity) cannot be omitted.

\medskip

\begin{prob}[Stable constraint maps]
    For suitable targets, such as when $M=\R^m\setminus \mathbb B^m$ for $m\geq 4$, is there a suitable notion of stability for constraint maps so that properties (i)--(iv) of Definition \ref{def:F} hold for such maps?
\end{prob}

This problem is closely related to the discussion in Section \ref{sec:prelims}.

\begin{prob}[Fine structure of FB]
    Consider the free boundary of a constraint map $u$ in the neighborhood of a continuity point $x_0\not\in \Sigma(u)$. Suppose also that $|Du(x_0)|>0$. Does the fine structure analysis of the free boundary from the scalar obstacle problem \cite{Figalli2019,FRS20,FZ25} carry through to the free boundary of $u$ near $x_0$?
\end{prob}

\medskip

In the next two problems, we consider specific examples of constraint maps and we inquire about their properties.

\begin{prob}[Singular FB produced by convex obstacles]
In \cite[\S 6]{FGKS2} we showed that for a large class of convex obstacles, and for a large class of boundary data, there are necessarily singularities of the corresponding minimizing constraint maps on their free boundaries. Is it the case that, at least in some of these examples, also the free boundaries are singular?
\end{prob}

\begin{prob}[Minimality of the vortex]
Is it the case that the radial map considered in Proposition \ref{prop:radial} and in \cite[\S 2.4]{FGKS2} is a minimizing constraint map in low dimensions, specifically $3\leq n\leq 6$?
\end{prob}

\medskip

Our final problem pertains to the more geometric setting of constraint maps between Riemannian manifolds, as described in \cite{D}. Here and in \cite{FKS,FGKS2,FGKS3} we only considered the case where the target manifold is flat and all of the non-trivial geometry is present on its boundary. Perhaps the simplest setting where also $\textup{int}(M)$ has non-trivial geometry is described in the following problem:

\begin{prob}[Constraint maps into a spherical cap]\label{prob:cap}
    For $h\in (-1,0]$, consider the spherical cap 
    $$\mb S^{m-1}_{h,+}\equiv\{x\in \mb S^{m-1}: x_n\geq h\}\subset \R^m
    $$
    whose (relative) boundary is the latitude line $P_h:=\partial \mb S^{m-1}_{h,+}$.
    We can consider energy minimizing maps $u\in W^{1,2}(\mb B^n, \mb S^{m-1}_{h,+})$, which develop a free boundary $\partial u^{-1}(P_h).$ Are there discontinuities of $u$ on the free boundary?
\end{prob}

We note that the answer to Problem \ref{prob:cap} may depend on $n$ and $h$: if $h=0$ and $n\leq 7$ then $u$ is smooth \cite{Schoen1984}. Moreover, the case $h\geq 0$ is likely not very interesting, since if $u$ omits a neighborhood of the equator it should be smooth \cite{Solomon1985}.

\section{Acknowledgements}
A. F. is grateful to the Marvin V. and Beverly J. Mielke Fund for supporting his stay at IAS Princeton, where part of this work has been done.
A.G. was partially supported by Dr.~Max Rössler, the Walter Haefner Foundation, and the ETH Zürich Foundation, as well as by a Newton International Fellowship from the Royal Society. 
S.K. was supported by the postdoctoral fellowship from the Verg Foundation and the starting grant from the Swedish Research Council (grant no.~2024-04747).
H.S. was supported by the Swedish Research Council (grant no.~2021-03700).

\section*{Declarations}

\noindent {\bf  Data availability statement:} All data needed are contained in the manuscript.

\medskip

\noindent {\bf  Funding and/or Conflicts of interests/Competing interests:} The authors declare that there are no financial, competing or conflict of interests.


\bibliographystyle{abbrv}
\bibliography{library.bib}

\end{document}